%Authors: P.G. Casazza and N.J. Kalton

%Title: Generalizing the Paley-Wiener perturbation theory for Banach spaces

%Filename: casazzakaltonpaley.tex
%TeX: AMSTeX
%Length: 26704 bytes
%Received Date: 9/8/97
%SubjectClass: 46B99
%Abstract: We extend the Paley-Wiener pertubation theory to linear operators
%mapping a subspace of one Banach space into another Banach space.

%Citation: Proc. Amer. Math. Soc.

%Special character check block
%32   space        33 ! exclam. pt.   34 " double quote  35 # sharp
%36 $ dollar       37 % percent       38 & ampersand     39 ' prime
%40 ( left paren.  41 ) rt. paren.    42 * asterisk      43 + plus
%44 , comma        45 - minus         46 . period        47 / divide
%58 : colon        59 ; semi-colon    60 < less than     61 = equal
%62 > greater than 63 ? question mark 64 @ at
%91 [ left bracket 92 \ backslash     93 ] right bracket 94 ^ caret
%95 _ underline    96 ` left single quote
%123 { left brace  124 | vertical bar 125 } right brace  126 ~ tilda

%Insert your TeX file starting here.

%This is AMS-TeX 2.1
%
\documentstyle{amsppt}
\magnification=\magstep1
\NoRunningHeads
\NoBlackBoxes
\parindent=1em
\vsize=7.4in
%macros
%
\topmatter

\title
Generalizing the Paley-Wiener Perturbation Theory for Banach Spaces
\endtitle
\author
Peter G. Casazza and Nigel J. Kalton
\endauthor
\address
Department of Mathematics,
The University of Missouri,
Columbia, Missouri 65211,
USA
\endaddress

\email
pete\@casazza.math.missouri.edu;   nigel\@math.missouri.edu
\endemail

\thanks
The first author was supported by NSF-DMS 9201357, the Danish
Natural Science Research Council, grant no. 9401598, and grants from
the University of Missouri Research Board, and the University of
Missouri Research Council.  The second author was supported by
NSF-DMS 95000125.
Paley-Wiener Perturbation Theory, spectrum, approximate fixed
points.
\endthanks

%\subjclass
%\endsubjclass
%
\abstract
We extend the Paley-Wiener pertubation theory to linear operators
mapping a subspace of one Banach space into another Banach space.
\endabstract

\endtopmatter
\document
\baselineskip=15pt
\heading{1. Introduction}
\endheading
\vskip10pt

In his classical book on potential theory in 1877, Carl Von Neumann
\cite{13} introduced what we now call the {\bf Neumann series} for
a linear operator:   If $X$ is a Banach space and $T:X\rightarrow X$
is a linear operator satisfying $\|I - T\| < 1$, then $T$ is an
onto isomorphism and $T^{-1} = \sum_{n=0}^{\infty}(I-T)^{n}$.
Special cases of this result were rediscovered by Paley and Wiener
in 1934 \cite{10} and in 1940 by Boas \cite{3}.  After further
generalizations by Pollard \cite{11} and Sz. Nagy \cite{9}, Hilding
\cite{6} gave the most general form:  If $X$ is a Banach space, and
$T:X\rightarrow X$ is a linear operator, $\lambda \in [0,1)$, and
for all $x\in X$, $\|(I - T)x\| \le {\lambda}(\|x\| + \|Tx\|),$ then
$T$ is an onto isomorphism.  We will investigate the more general
setting where $Y$ is a subspace of a Banach space $X$, $Z$ is a
Banach space and $S,T:Y\rightarrow Z$ are linear operators
satisfying for all $x\in X$, the inequality, $\|Sx - Tx\| \le
{{\lambda}_{1}}\|Sx\| + {{\lambda}_{2}}\|Tx\|$, where
${{\lambda}_{1}}, {{\lambda}_{2}} \in [0,1)$.  In this case,
properties of $S$ will carry over to $T$.
This includes being one-to-one, onto, closed, open, having dense
range, being a quotient map and most importantly, being a Fredholm
operator (and the Fredholm index is maintained).  A special case of
this result is a generalization of the theorem of Neumann:  If $Y$
is a subspace of a Banach space $X$ and $T:Y\rightarrow X$ is a
linear operator with $\|(I - T)|_{Y}\| < 1$, then $Y$ and $TY$ have
the same codimension in $X$.

\heading{2. The Basic Inequalities}
\endheading
\vskip10pt

We will first develop the basic inequalities needed throughout the
paper.
We will always assume that ${{\lambda}_{1}}, {{\lambda}_{2}}$ are
real numbers with ${{\lambda}_{1}}, {{\lambda}_{2}} \in [0,1)$.

\proclaim{Basic Inequality}
If $x,y$ are elements of a Banach space $X$ satisfying,
$$
\|x - y\| \le {{\lambda}_{1}}\|x\| +  {{\lambda}_{2}}\|y\|,  \tag 1
$$
then,
$$
 \frac{1 - {{\lambda}_{2}}}{1 + {{\lambda}_{1}}}\|y\| \le \|x\| \le
\frac{1 + {{\lambda}_{2}}}{1 - {{\lambda}_{1}}}\|y\|. \tag 2
$$
\endproclaim

\demo{Proof}
With $x,y$ as above,
$$
\|x\| \le \|x - y\| + \|y\| \le {{\lambda}_{1}}\|x\| +
{{\lambda}_{2}}\|y\| + \|y\|.
$$
It follows that,
$$
\|x\| \le \frac{1 + {{\lambda}_{2}}}{1 - {{\lambda}_{1}}}\|y\|.
$$
Switching the roles of $x$ and $y$ above completes the inequality.
\enddemo

We will be working in the case where $Y,Z$ are subspaces of a
Banach space $X$ and $T:Y\rightarrow Z$ is a bounded linear
operator.  Letting $y = Tx$ in $(1)$, we have for all $x\in Y$,
$$
\|(I - T)x\| \le {{\lambda}_{1}}\|x\| + {{\lambda}_{2}}\|Tx\|. \tag 3
$$

It follows from our basic inequality that $T$ is an isomorphism of
$Y$ into $X$.  We now show that if $T$ satisfies $(3)$, then so do
certain operators obtained from $T$.

\proclaim{Proposition 1}
Let $Y,Z$ be subspaces of a Banach space $X$ and $T:Y\rightarrow Z$
be a surjective linear operator satisfying for all $x\in Y$,
$$
\|(I - T)x\| \le {{\lambda}_{1}}\|x\| + {{\lambda}_{2}}\|Tx\|,
$$

Then

(1) $T^{-1}$ satisfies for all $x\in Z$,
$$
\|(I - T^{-1})x\| \le {{\lambda}_{2}}\|x\| +
{{\lambda}_{1}}\|T^{-1}x\|. \tag 4
$$

(2)  For every $\alpha > 0$, ${\alpha}T$ satisfies for all $x\in Y$,
$$
\|(I - {\alpha}T)x\| \le {{\lambda}^{'}_{1}}\|x\| +
{{\lambda}^{'}_{2}}\|{\alpha}Tx\|, \tag 5
$$
 with constants
${{\lambda}^{'}_{1}} = \text{max}\{1 - {\alpha}(1 -
{{\lambda}_{1}}, {{\lambda}_{1}}\}$, and ${{\lambda}^{'}_{2}} =
\text{max}\{1 - \frac{1 + {{\lambda}_{2}}}{\alpha},
{{\lambda}_{2}}\}$.

(3)  For every $0 \le \alpha \le 1$, the operator  $T_{\alpha} = (1
-{\alpha})I + {\alpha}T$ satisfies for all $x\in Y$,
$$
\|(I - T_{\alpha})x\| \le {{\lambda}^{'}_{1}}\|x\| +
{{\lambda}^{'}_{2}}\|T_{\alpha}x\|,  \tag 6
$$
where ${{\lambda}^{'}_{1}} = {\alpha}{{\lambda}_{1}} + (1 -
{\alpha}){{\lambda}_{2}}$, and ${{\lambda}^{'}_{2}} =
{\lambda}_{2}$.

(4)  For every $a < 0$, the operator $aI-T$ is an isomorphism.
\endproclaim

\demo{Proof}
(1)  If $z = Tx$, then
$$
\|(I - T^{-1})z\|_{Z\rightarrow Y} = \|Tx - x\| \le
{{\lambda}_{1}}\|x\| + {{\lambda}_{2}}\|Tx\| =
{{\lambda}_{1}}\|T^{-1}z\| + {{\lambda}_{2}}\|z\|.
$$

(2)  Although this can be done in one case, we will do it in two
cases to identify the exact constants obtained in each case.
\proclaim{Case I}
Assume $\alpha \le 1$.
\endproclaim

For any $x\in Y$,
$$
\align
\|(I - {\alpha}T)x\| &= \|(1 - {\alpha})I + {\alpha}(I - T)y\| \\
 &\le (1 - {\alpha})\|x\| + {{\lambda}_{1}}{\alpha}\|x\| +
{{\lambda}_{2}}{\alpha}\|Tx\| \\
&= [1 - {\alpha}(1 - {{\lambda}_{1}})]\|x\| +
{{\lambda}_{2}}\|{\alpha}Tx\|.
\endalign
$$

\proclaim{Case II}
Assume $\alpha > 1$.
\endproclaim

For any $x\in Y$,
$$
\align
\|x - {\alpha}Tx\| &\le \|(I - T)x\| + (\alpha - 1)\|Tx\| \\
&\le {{\lambda}_{1}}\|x\| + {{\lambda}_{2}}\|Tx\| + (\alpha -
1)\|Tx\| = {{\lambda}_{1}}\|x\| + \frac{{\lambda}_{2} + \alpha -
1}{\alpha}\|{\alpha}Tx\|.
\endalign
$$

(3)  For any $x\in Y$,
$$
\align
\|(I - T_{\alpha})x\| &= {\alpha}\|(I - T)x\| \le
{\alpha}{{\lambda}_{1}}\|x\| + {\alpha}{{\lambda}_{2}}\|Tx\|
= {\alpha}{{\lambda}_{1}}\|x\| + {{\lambda}_{2}}\|{\alpha}Tx\| \\
&\le {\alpha}{{\lambda}_{1}}\|x\| + {{\lambda}_{2}}\|{\alpha}Tx +
(1 - {\alpha})x\| + {{\lambda}_{2}}(1 - {\alpha})\|x\| \\
&= [{\alpha}{{\lambda}_{1}} + (1 - {\alpha}){{\lambda}_{2}}]\|x\| +
{{\lambda}_{2}}\|T_{\alpha}x\|.
\endalign
$$

(4)  This is immediate from (3) and the observation,
$$
aI-T = -(|a|+1)[(1-\frac{1}{|a|+1})I + \frac{1}{|a|+1}T].
$$

\enddemo

If we weaken inequality 3, for example by letting ${{\lambda}_{1}}
= 1$, we lose our conclusion that $T$ is a bounded linear operator.
For example, it is immediate that $T:{\ell_p}\rightarrow {\ell_p}$
given by $T(\{a_{i}\}) = \{ia_{i}\}$, satisfies $\|(I - T)x\| \le
\|Tx\|$, for all $x\in {\ell}_{p}$.  Also, $T = 0$ satisfies
inequality 3 with ${{\lambda}_{1}} = 1$.  The next proposition shows
that this is essentially all that can go wrong with the weaker
inequality 3.

\proclaim{Proposition 2}
Let $Y$ be a subspace of a Banach space $X$ and $T:Y\rightarrow X$
be a linear map.

(1)  Suppose $T$ is bounded and there exists ${\lambda}\in [0,1)$
such that for all $x\in Y$,
$$
\|x - Tx\| \le {\lambda}\|x\| + \|Tx\|.
$$
Choose ${\delta}\in [{\lambda},1)$ so that $\|T\| \le \frac{\delta
- {{\lambda}}}{1 - \delta}$.  Then $T$ satisfies for all $x\in Y$,
$$
\|x - Tx\| \le {\delta}(\|x\| + \|Tx\|).
$$

(2)  If $T^{-1}$ is bounded on $T(Y)$, and $T$ satisfies for all
$x\in Y$,
$$
\|x - Tx\| \le \|x\| + {{\lambda}_{2}}\|Tx\|,
$$
then
$$
\|x - Tx\| \le {\delta}(\|x\| + \|Tx\|)
$$
where ${\delta}\in [{\lambda},1)$ is chosen so that  $\|T^{-1}\|
\ge \frac{1 - {\delta}}{\delta - {{\lambda}}} > 0.$
\endproclaim

\demo{Proof}
(1)  Since $\|Tx\| \le \frac{\delta - {{\lambda}}}{1 -
\delta}\|x\|$, we see that $(1 - {\delta})\|Tx\| \le ({\delta} -
{{\lambda}})\|x\|$.  Hence,
${{\lambda}}\|x\| + \|Tx\| \le {\delta}(\|x\| + \|Tx\|)$.  (2)
follows similarly.
\enddemo

\heading{3.  An Application from Operator Theory}
\endheading
\vskip10pt

We denote the unit sphere of a Banach space $X$ by $S_{X} = \{x\in
X: \|x\| = 1\}$.  Let ${\sigma}(T)$ denote the {\bf spectrum} of an
operator $T:X \rightarrow X$ and ${\pi}(T)$ denote the {\bf
approximate point spectrum} of $T$.  That is, $\lambda \in
{\sigma}(T)$ if $T-{\lambda}I$ is not invertible, and ${\lambda} \in
{\pi}(T)$ if there is a sequence $x_{n} \in S_{X}$ so that
$\|(T-{\lambda}I)x_n\| \rightarrow 0$.  In this terminology,
proposition 7.9 in \cite{7} states,

\proclaim{Theorem (Kalton, Peck, Roberts)}
The complement of the spectrum of $T$ is a clopen (i.e. both closed
and open) set in the complement of the approximate point spectrum
of $T$, which contains the unbounded component.
\endproclaim

\proclaim{Corollary 3}
Let $X$ be a Banach space and $T:X\rightarrow X$ be an isomorphism
of $X$ into $X$.

(1)  If the operator ${\alpha}I - T$ is an isomorphism for all
${\alpha} > 0$, then $T$ is onto.

(2)  If the operator ${\alpha}I - T$ is an isomorphism for all
${\alpha} < 0$, then $T$ is onto.
\endproclaim

\demo{Proof}
By our assumption in $(1)$, $[0,\infty)$ is a subset of the
unbounded component of ${{\pi}(T)}^{c}$ and hence is a subset of
${{\sigma}(T)}^{c}$.  (2) follows similarly.
\enddemo

Some assumption on $T$ in Corollary 3 is necessary, since without
it $T = 0$ satisfies the hypotheses. In the complex case, the
hypotheses in theorem 3 could be stated more generally as: (1)  If
for some complex unit $|{{\lambda}_{0}}| = 1$, we have that
${\alpha}{{\lambda}_{0}}I - T$ is an isomorphism for all ${\alpha} >
0$, then $T$ is onto.  Similarly for (2).

The classical Borsuk-Ulam theorem asserts that any continuous map
from an n-dimensional sphere to itself must either be onto, or have
both fixed points and antipodal points.  The Borsuk-Ulam theorem
fails for infinite dimensional Banach spaces, in its exact form,
even for linear isometries (just take the shift operator on a
Hilbert space).  However, there is an approximate version of this
theorem.  We say that a sequence of elements $\{x_n\}$ in $S_{X}$ is
an {\bf approximate fixed point sequence} for a mapping
$f:S_{X}\rightarrow X$ if $\lim_{n\rightarrow \infty}\|x_n -
f(x_{n})\| = 0$.  It is an {\bf approximate antipodal sequence} for
$f$ if $\lim_{n\rightarrow \infty}\|-x_n - f(x_{n})\| = 0$.
Benyamini and Sternfeld \cite{2} have shown that every infinite
dimensional Banach space $X$ has a Lipschitz map of the unit ball of
$X$ into itself without approximate fixed points.  If $X = \ell_p$,
for $1 < p < \infty$, this map automatically satisfies inequality
3.
But the theorem of Kalton, Peck, and Roberts above does yield an
approximate version of the Borsuk-Ulam theorem for linear
isometries.

\proclaim{Borsuk-Ulam Theorem for Linear Operators}
If $X$ is a Banach space and $T:X\rightarrow X$ is an isometry,
then either $T$ is onto, or $T$ has both an approximate fixed point
sequence and an approximate antipodal sequence.  Moreover, if $X$ is
complex and $T$ is an isometry which is not onto, then the spectrum
of $T$ contains the unit circle.
\endproclaim

\demo{Proof}
All of this is immediate from Corollary 3 except the last statement
which only requires the observation that if $T$ is an isometry
which is not onto, then for all complex numbers $|{\lambda}| = 1$,
the operator ${\lambda}T$ is not onto.  By Corollary 3, there is an
${\alpha} > 0$ so that ${\alpha}I - {\lambda}T$ is not an
isomorphism.  It follows that ${\alpha} = 1$, so $I - {\lambda}T$ is
not an isomorphism, and hence $\overline{\lambda}I - T$ is not an
isomorphism.
\enddemo

Isometries which are onto need not have approximate fixed points or
approximate antipodal points.  To see this, define $T:\ell_p
\rightarrow \ell_p$ by:
$$
T(a_{1}, a_{2}, a_{3},a_{4},\ldots ) =
(a_{2},-a_{1},a_{4},-a_{3},\ldots).
$$
Then $T$ is an isometry of $\ell_p $ onto itself, but $T$ has no
approximate fixed points or approximate antipodal points.  Hence,
${\alpha}I - T$ is an onto isomorphism for all ${\alpha}$.

It is clear that any operator satisfying inequality (3) cannot have
approximate fixed points or approximate antipodal points.  One
would hope that an operator which fails inequality (3) would need to
have approximate fixed points or approximate antipodal points.  It
is easily checked that this is the case in a uniformly convex space.
 To see that this is not true in general, let $T(a_{1},a_{2}) =
(a_{2},-a_{1})$ be considered as an operator on $\ell^2_1$.

The classical perturbation theorem of Hilding \cite{6} now follows.

\proclaim{Hilding's Pertubation Theorem}
If $X$ is a Banach space and $T:X\rightarrow X$ satisfies, for all
$x\in X$,
$$
\|(I - T)x\| \le {{\lambda}_{1}}\|x\| + {{\lambda}_{2}}\|Tx\|,
$$
for some ${{\lambda}_{1}}, {{\lambda}_{2}} \in [0,1)$, then $T$ is onto.
\endproclaim

\demo{Proof}
Proposition 1 (4) states that $T$ satisfies hypothesis 2 of Corollary 3.
\enddemo

\heading{4. Generalizing Paley-Wiener}
\endheading
\vskip10pt

Now we will extend the theory to operators between subspaces of
Banach spaces.  We start with an elementary observation.

\proclaim{Lemma 4}
Let $X,Y$ be Banach spaces and $S,T:X\rightarrow Y$ be linear
operators satisfying,
$$
\|Sx - Tx\| \le {{\lambda}_{1}}\|Sx\| + {{\lambda}_{2}}\|Tx\|,
$$
for all $x\in X$, and  fixed ${{\lambda}_{1}}, {{\lambda}_{2}} \in
[0,1)$.  Then if $S$ has closed range (respectively, is one-to-one,
has dense range, is an open map, is a quotient map, is an
isomorphism) then $T$ has closed range (respectively, is one-to-one,
has dense range, is an open map, is a quotient map, is an
isomorphism).
\endproclaim

\demo{Proof}
Applying our Basic Inequality to $Sx,Tx$ we have:
$$
\frac{1 - {{\lambda}_{2}}}{1 + {{\lambda}_{1}}}\|Tx\| \le \|Sx\|
\le \frac{1 + {{\lambda}_{2}}}{1 - {{\lambda}_{1}}}\|Tx\|.
$$
It follows that $L:S(X)\rightarrow T(X)$ defined by, $L(Sx) = Tx$,
is a well defined onto isomorphism, which therefore has a unique
extension to an isomorphism of $\overline{S(X)}$ onto
$\overline{T(X)}$.  This is all that is needed for the proof of the
theorem.
\enddemo

Now we want to prove a deeper generalization of the Paley-Wiener
perturbation theory.  We will need a result of Krein, Krasnoselskii,
and Milman \cite{8},  which can also be found in Gokhberg and Krein
\cite{5} or Day \cite{4}.

\proclaim{Lemma 5}
Let $E,F$ be subspaces of a Banach space $X$ with $\text{dim } F <
\infty$ and $\text{dim }E > \text{dim } F$.  Then there exists an $0
\not= x \in S_{E}$ such that
$$
1 = \|x\| = d(x,F).
$$
\endproclaim

To prove our main result, we need two lemmas.  The first is
actually a special case of the result.

\proclaim{Lemma 6}
Let $Y$ be a subspace of a Banach space $X$, and $T:Y\rightarrow X$
a linear operator satisfying $\|(I - T)|_{Y}\| < 1$.  Then for
every subspace $W$ in $Y$, $\text{codim}_{X}W \ge
\text{codim}_{X}T(W).$
\endproclaim

\demo{Proof}
Let $\lambda = \|(I - T)|_{Y}\|.$  If the lemma fails, then there
is a subspace $W$ in $Y$ so that $\text{dim}(TW)^{\perp} >
\text{dim}W^{\perp}$.  By lemma 5, there is an element $x^{*}\in
S_{(TW)^{\perp}}$ with $d(x^{*},W^{\perp}) = 1 =
\text{sup}\{x^{*}(x):x\in B_{W}\}.$  Now, for any $\delta > 0$,
there is an $x\in B_{W}$ so that $x^{*}(x) \ge 1 - {\delta}$.  Since
$x^{*}(T(W)) = 0$, we have,
$$
1 - {\delta} \le x^{*}(x) = x^{*}(x - Tx) \le
{\lambda}\|x^{*}\|\|x\| = {\lambda}.
$$
But, this is a contradiction for $1 - \delta > {\lambda}$.  Therefore,
$$
\text{codim}_{X}W = \text{dim}W^{\perp} \ge \text{dim}(TW)^{\perp}
= \text{codim}_{X}T(W).
$$
\enddemo

\proclaim{Lemma 7}
Let $Y$ be a subspace of a Banach space $X$ and $T:Y\rightarrow X$
be a linear operator with $\|(I - T)|_{Y}\| < \frac{1}{2}$.  Then,
$\|(I - T^{-1})|_{TY}\| < 1$.
\endproclaim

\demo{Proof}
For any $x\in X$,
$$
\|Tx\| \ge \|x\| - \|(I - T)x\| \ge \frac{1}{2}\|x\|.
$$
It follows that $\|T^{-1}|_{TY}\| \le 2.$ Now, for every $y\in TY$,
$$
\|(I - T^{-1})y\| = \|(I - T)T^{-1}y\| \le \|I - T\|\|T^{-1}y\| \le
\|I - T\| 2\|y\|.
$$
\enddemo

We can now prove the main result of this section.  Again, the proof
is inspired by Hilding \cite{6}.

\proclaim{Theorem 8}
Let $Y$ be a subspace of a Banach space $X$,
${\lambda}_{1},{\lambda}_{2}\in [0,1)$, and $T:Y\rightarrow X$ a
linear operator satisfying, for all $x\in Y$,
$$
\|(I - T)x\| \le {{\lambda}_{1}}\|x\| + {{\lambda}_{2}}\|Tx\|.
$$
Then $\text{codim}_{X}Y = \text{codim}_{X}T(Y)$.
\endproclaim

\demo{Proof}
To simplify the proof, let ${\lambda} =
\text{max}\{{{\lambda}_{1}},{{\lambda}_{2}}\}$.  With $T_{\alpha}$
defined as in Proposition 1 (3), we have from our Basic Inequality,
for all $0 \le {\alpha} \le 1$,
$$
\frac{1 - {\lambda}}{1 + {\lambda}} \|x\| \le \|T_{\alpha}x\|. \tag 7
$$
Next, we let
$$
E = \{0 \le {\alpha} \le 1 : \text{codim}_{X}T_{\alpha}Y =
\text{codim}_{X}Y\}.
$$

If $\alpha = 0$, then $T_{\alpha} = I$, so $0\in E \not= \phi$.

Next, we will show that for all ${\alpha}$ sufficiently close to
${\beta}$, we have $\text{codim}_{X}T_{\alpha}Y =
\text{codim}_{X}T_{\beta}Y$.
To see this, given $0 \le {\alpha},{\beta} \le 1$, we have
$$
\|T_{\alpha}x - T_{\beta}x\| = \|({\alpha} - {\beta})Tx - ({\alpha}
- {\beta})x\| \le |{\alpha} - {\beta}|(1 + \|T\|)\|x\|.
$$

Hence if we let
$$
\epsilon = \frac{1}{2} \frac{1 - \lambda}{1 + \lambda} \frac{1}{1 +
\|T\|},
$$
then as long as $|{\alpha} - {\beta}| \le {\epsilon}$, and applying
(7) we have
$$
\|T_{\alpha}x - T_{\beta}x\| \le
\frac{1}{2\|T^{-1}_{\alpha}\|}\|x\|.  \tag 8
$$
Now define an operator $L:T_{\alpha}(X) \rightarrow T_{\beta}(X)$
by $LT_{\alpha}x = T_{\beta}x$.  By inequality (8), we have that
$\|I - L\| < \frac{1}{2}$.   Thus by lemma 7, we can apply lemma 6
to both $L$ and $L^{-1}$ to conclude that
$\text{codim}_{X}T_{\alpha}Y = \text{codim}_{X}T_{\beta}Y$.

Summarizing, we have that $0\in E$, and whenever ${\alpha} \in E$,
we have that
$({\alpha} - {\epsilon},{\alpha} + {\epsilon}) \cap [0,1] \subset
E$.  Hence, $E = [0,1]$ and so $1 \in E$, which is what we needed.
\enddemo

Theorem 8 gives a generalization of the result of Neumann \cite{13}.

\proclaim{Corollary 9}
If $Y,Z$ are subspaces of a Banach space $X$, and if
$T:Y\rightarrow Z$ is a surjective linear operator with $\|(I -
T)|_{Y}\| < 1$, then $\text{codim}_{X}Y = \text{codim}_{X}Z$.
\endproclaim

We could obtain Corollary 9 directly from the Neumann series if $T$
had an extension $\hat{T}$ to all of $X$ also satisfying $\|I -
\hat{T}\| < 1$.  In fact, we can get such an extension if there is a
projection $P:X\rightarrow Y$ with $\|P\| < 1$ (or just $\|P\| <
\|I - T\|$).  In this case we define $\hat{T}:X\rightarrow X$ by
$\hat{T}x = TPx + (I-P)x$.  Now, for all $x\in X$,
$$
\|(I - \hat{T})x\| = \|(I - T)Px\| \le \|I - T\|\|Px\| \le \|I -
T\|\|x\|.
$$
However, in general $T$ need not have an extension which is an
isomorphism on $X$. Our next example shows that even if $\text{dim
}Y < \infty$, there need not be an extension $\hat{T}$ of $T$
satisfying $\|I - \hat{T}\| < 1$.

\proclaim{Example 10}
Let $X = {\ell_p}\oplus {\ell_p}$, for $p \not= 2$.  Choose a
subspace $W$ in $\ell_p$ which is isomorphic to $\ell_p$ but
uncomplemented in $\ell_p$.  Let $Y = W \oplus 0$, and let
$\{({f_n},0)\}$ be the unit vector basis of $\ell_p$ in $Y$, and
$\{e_n\}$ be the unit vector basis of $\ell_p$. Also, let
$$
Z = \{(\sum a_{n}f_{n},\frac{1}{K} {\sum a_{n}e_{n}) : (a_n)\in
{\ell_p}} \},
$$
where $K$ is chosen so that,
$$
\frac{1}{K} \|\sum a_{n}e_{n}\| \le \frac{1}{2} \|\sum a_{n}f_{n}\|.
$$
Finally, define $T:Y\rightarrow Z$ by
$$
T(\sum a_{n}f_{n},0) = (\sum a_{n} f_{n}, \frac{1}{K} \sum
a_{n}e_{n}). \tag 9
$$
Then,

(1)  $\|I - T\| < 1$,

(2)  $Z$ is complemented in $X$,

(3)  $T$ is an isomorphism of $Y$ onto $Z$.

Therefore, $T$ cannot be extended to be an isomorphism of $X$ onto $X$.
\endproclaim

\demo{Proof}
(1) We just apply (9) to get
$$
\|(I - T)(\sum a_{n} (f_{n},0))\|_{X} = \frac{1}{K} \|\sum
a_{n}e_{n}\|_{\ell_p}   \le \frac{1}{2}\|\sum a_{n}f_{n}\|.
$$

(2) we define the operator $P$ on $X$ by:
$$
P(\sum a_{n}e_{n}, \sum b_{n}e_{n}) =  (K \sum b_{n}f_{n},\sum
b_{n}e_{n}).
$$
It is clear that $P$ is a bounded linear projection of $X$ onto $Z$.

(3) This is clear since the operator $T(f_{n}) = (f_{n},e_{n})$ is
an isomorphism.
\enddemo

By a standard compactness arguement, we can finite dimensionalize
the above example.  There is a choice of natural numbers $j_{1} <
j_{2} < j_{3} < \cdots$ with the following property.  Let $Y_{n} =
\text{span}_{1 \le i \le j_n}(f_{i},0)$ and $Z_{n} = \text{span}_{1
\le i \le j_n}(f_{i},\frac{1}{K}e_{i})$, and $T_{n} = T|_{Y_n}$.
Then there is a $0 < \lambda < 1$ so that for all $n =
1,2,3,\ldots$, $\|I - T_{n}\| < {\lambda}$, but for any extension
$\hat{T_n}$ of $T_{n}$ to all of $X$, we have
$\|\hat{T}_{n}\|\|\hat{T}^{-1}_{n}\| \ge n$.

The next corollary of Theorem 8 comes from mimmicking the proof of
lemma 4.

\proclaim{Corollary 11}
If $X,Y$ are Banach spaces and $S,T:X\rightarrow Y$ are linear
operators satisfying for all $x\in X$,
$$
\|Sx - Tx\| \le {\lambda}(\|Sx\| + \|Tx\|),
$$
then $\text{codim}_{Y}(\overline{S(X)}) =
\text{codim}_{Y}(\overline{T(X)})$.
\endproclaim

Given an operator $T:X\rightarrow Y$ with closed range, we let
$$
{\alpha}(T) = \text{dim }\text{ker } T
$$
and
$$
{\beta}(T) = \text{dim }Y/T(X).
$$
If either ${\alpha}(T) < {\infty}$ or ${\beta}(T) < {\infty}$ we
define the {\bf Fredholm index i(T)} of $T$ by $i(T) = {\alpha}(T) -
{\beta}(T)$.  If ${\alpha}(T)$ and ${\beta}(T)$ are both finite
(i.e. if $i(T)$ is defined and is finite) then $T$ is called a {\bf
Fredholm operator} of index $i(T)$.

\proclaim{Corollary 12}
Let $X,Y$ be Banach spaces and $S,T:X\rightarrow Y$ be linear
operators satisfying for all $x\in X$, and fixed $0 \le
{{\lambda}_{1}},{{\lambda}_{2}} < 1$,
$$
\|(S - T)x\| \le {{\lambda}_{1}}\|Sx\| + {{\lambda}_{2}}\|Tx\|.
$$
If $S$ is a Fredholm operator with Fredholm Index $n$, then $T$ is
also a Fredholm operator with Fredholm Index $n$.
\endproclaim

\demo{Proof}
By our Basic Inequality, $\text{ker }S = \text{ker }T$.  Now apply
Theorem 8.
\enddemo

We end with one final application of Theorem 8.

\proclaim{Corollary 13}
Let $X$ be a Banach space and $T:X\rightarrow X$ be a linear
operator satisfying
$$
\|(I - T)x\| \le {\lambda}(\|x\| + \|Tx\|),  \tag 10
$$
for all $x\in X$, and fixed $0 \le {\lambda} < 1$.  Then, for all
natural numbers n, and all $x\in X$, we have, $x\in
\overline{\text{span}_{n\le k}T^{k}x}.$
\endproclaim

\demo{Proof}
For each $n = 0,1,2,\cdots$, let $W_{n} =
\overline{\text{span}_{n\le k}T^{k}x}.$  Then $T$ maps $W_{n}$ into
$W_{n}$ and satisfies (10), and hence is onto.  Since $T(W_{n}) =
W_{n+1}$, we see that $W_{0} = W_{1} = W_{2} = \cdots$, which proves
the corollary.
\enddemo

\proclaim{ACKNOWLEDGEMENT}
The authors thank Ole Christensen for helpful discussions
concerning this paper.
\endproclaim


\Refs

\ref\no{1}
\by M.G. Arsove
\paper The Paley-Wiener theorem in metric linear spaces
\jour Pacific J. Math
\vol 10
\yr 1960
\pages 365-379
\endref

\ref\no{2}
\by Y. Benyamini and Y. Sternfeld
\paper Spheres in infinite-dimensional normed spaces are Lipschitz
contractible
\jour Proceedings of the American Mathematical Society
\vol 88 no. 3
\yr 1983
\pages 439-445
\endref


\ref\no{3}
\by R.P. Boas, Jr.
\paper Expansions of analytic functions
\jour Transactions of the American Mathematical Society
\vol 48
\yr 1940
\pages 467-487
\endref

\ref\no{4}
\by M.M. Day
\paper On the basis problem in normed spaces
\jour Proceedings of the American Mathematical Society
\vol 13
\yr 1962
\pages 655-658
\endref

\ref\no{5}
\by I.C. Gokhberg and M.G. Krein
\paper Fundamental theorems on deficiency numbers, root numbers,
and indices of linear operators (Russian)
\jour Usp. Mat. Nauk
\vol 12
\yr 1957.
\pages 43-118.  Translated in American Mathematical Society
Translations, ser. 2,
vol. 13
\endref


\ref\no{6}
\by S.H. Hilding
\paper Note on completeness theorems of Paley-Wiener type
\jour Annals of Math
\vol 49
\yr 1948
\pages 953-955
\endref

\ref\no{7}
\by N.J. Kalton, N.T. Peck, and J.W. Roberts
\paper An F-space sampler
\jour London Mathematical Society Lecture Notes, Cambridge
University Press
\vol 89
\yr 1984
\endref


\ref\no{8}
\by M.G. Krein, M.A. Krasnoselskii, and D.P. Milman
\paper On the defect numbers of linear operators in a Banach space
and on some geometrical questions (Russian)
\jour Sb. Tr. Inst. Mat. Akad. Nauk Ukr. SSR
\vol 11
\yr 1948
\pages 97-112
\endref


\ref\no{9}
\by B. Sz. Nagy
\paper Expansion theorems of Paley-Wiener type
\jour Duke Math. Journal
\vol 14
\yr 1947
\pages 975-978
\endref

\ref\no{10}
\by R.E.A.C. Paley and N. Wiener
\paper Fourier transforms in the complex domain
\jour New York
\yr 1934
\endref

\ref\no{11}
\by H. Pollard
\paper Completeness theorems of Paley-Wiener type
\jour Annals of Math
\vol 45
\yr 1944
\pages 738-739
\endref


\ref\no{12}
\by James R. Retherford
\paper Basic sequences and the Paley-Wiener criterion
\jour Pacific Journal of Math
\vol 14 no. 3
\yr 1964
\pages 1019-1027
\endref

\ref\no{13}
\by C. Von Newmann
\paper Untersuchungen uber das logarithmische und Newtonsche Potential
\jour Teubner, Leipzig
\yr 1877
\endref


\endRefs
\enddocument